\documentclass[a4paper,12pt,reqno]{article}
\usepackage[T1]{fontenc}
\usepackage{authblk}
\usepackage{epsfig,graphicx,subfigure}
\usepackage{amsthm,amsmath,amssymb,amsfonts}
\usepackage{color,xcolor}
\usepackage{hyperref}
\usepackage[all]{xy}
\usepackage[top=40mm, bottom=40mm, left=40mm, right=35mm]{geometry}
\usepackage{tikz-cd}
\usepackage{hyperref}

\newtheorem{theorem}{Theorem}[section]

\newtheorem{proposition}[theorem]{Proposition}
\newtheorem{corollary}[theorem]{Corollary}
\theoremstyle{definition}
\newtheorem{definition}[theorem]{Definition}
\newtheorem{example}[theorem]{Example}

\numberwithin{equation}{section}

\begin{document}
\setcounter{page}{1}
\title{\vspace{-1.5cm}
\vspace{.5cm}
\vspace{.7cm}
{\large{\bf  On $g-$Fusion Frames Representations via  Linear Operators} }}
 \date{}
\author{{\small \vspace{-2mm} S. Jahedi$^{1}$, F. Javadi$^1$, M. J. Mehdipour$^{1}$\footnote{Corresponding author}  }}
\affil{\small{\vspace{-4mm}  $^1$ Department of Mathematics, Shiraz University of Technology, P. O. Box $71555-313$, Shiraz,  Iran.}}
\affil{\small{\vspace{-4mm} jahedi@sutech.ac.ir}}
\affil{\small{\vspace{-4mm} f.javadi@sutech.ac.ir}}
\affil{\small{\vspace{-4mm} mehdipour@sutech.ac.ir}}
\maketitle
\hrule
\begin{abstract}
\noindent
Let $\{\frak{M} _k \} _{ k \in \mathbb{Z}} $ be a sequence of closed subspaces of Hilbert space $H$, and let $\{\Theta_k\}_{k \in \mathbb{Z}}$  be a sequence of linear operators from  $H$ into $\frak{M}_k$, $k \in \mathbb{Z}$. In the definition of fusion frames, we replace the orthogonal projections on $\frak{M} _k$ by $\Theta_k$ and find a slight generalization of fusion frames. In the case where, $\Theta_k$  is  self-adjoint and $\Theta_k(\frak{M} _k)= \frak{M} _k$
for all $k \in \mathbb{Z}$, we show that if a
$g-$fusion frame $\{(\frak{M} _k, \Theta_k)\}_{k \in \mathbb{Z}}$ is represented via a linear operator $T$ on
$\hbox{span} \{\frak{M} _k\}_{ k \in \mathbb{Z}}$, then $T$ is bounded;  moreover, if $\{(\frak{M} _k, \Theta_k)\}_{k \in \mathbb{Z}}$ is  a tight $g-$fusion frame, then $T$ is not invertible.
We also study the perturbation and the stability of these fusion frames.
Finally, we give some examples to show the validity of the results.
 \end{abstract}

\\
\hrule
\vspace{0.5 cm}
\baselineskip=.55cm
\section{Introduction}

Let $\{f_k\}$ be a frame in a Hilbert space $H$ and $T$ is a linear operator (not necessarly bounded) on $\hbox{span}\{f_k\}$ such that $f_k=T^kf_0$. Christensen et.al [4] considered these types of frames and studied properties of them. Especially, they gave a necessary and sufficient condition for boundedness of $T$. Several authors investigated the results obtained in [4, 5, 6] for the various types of frame's extensions, where $T$ is a bounded linear operator; see for example [11, 12, 13] for $g-$frames; also [8, 10] for $K$-frames. The purpose of this paper is to investigate a slightly generalized fusion frames $\{(\frak{M} _k, \Theta_k)\}_{k \in \mathbb{Z}}$ that is represented via a linear operator (not necessarily bounded), where $\Theta_k$ is a linear operator. 

In this paper, we present a generalization of the concept of fusion frames and then investigate their representability via a linear operator. 
 In Section 3,  we assume that $\Theta_k$  is  self-adjoint and $\Theta_k(\frak{M} _k)= \frak{M} _k$
for all $k \in \mathbb{Z}$ and show that if a
$g-$fusion frame $\{(\frak{M} _k, \Theta_k)\}_{k \in \mathbb{Z}}$ is represented via a linear operator $T$ on
$\hbox{span} \{\frak{M} _k\}_{ k \in \mathbb{Z}}$, then $T$ is bounded;  also, if $\{(\frak{M} _k, \Theta_k)\}_{k \in \mathbb{Z}}$ is a  tight $g-$fusion frame, then $\{\Theta_k\}_{k \in \mathbb{Z}}$  is not represented via an invertible linear operator on  $\hbox{span}\{\frak{M} _k\}_{k \in \mathbb{Z}}$.  In Section 4, we show that the perturbation of representable fusion frames via a linear operator can also be represented by a linear operator under some suitable conditions. The stability of this type of fusion frame will be investigated. Finally, in Section 5, some examples are presented to support the validity of the main results. 
\section{Preliminaries}
Throughout this paper,  the space  $H$  denotes an  infinite dimensional separable Hilbert space,
$\{\frak{M} _k\}_{ k \in \mathbb{Z} }$  is a sequence of closed
 subspaces of $H$ and $\mathcal{L}(H)$ is  the space  of  all bounded linear operators on  $H$.
Let  $\Theta_0 $ be an  operator in  $\mathcal{L}(H)$ such that  $\Theta_0 (H) \subseteq \hbox{span} \{\frak{M} _k \}_{k \in \mathbb{Z}} $ and $\mathcal{T}$ be the right-shift operator on
\begin{eqnarray} \nonumber
\mathcal{H}= \bigg \{ \{f_k\}_{k\in \mathbb{Z}} \,: \,\,   f_k \in \frak{M} _k , \,\,\,\,\,  \sum_{k \in \mathbb{Z}} \| f_k\|^2 < \infty \bigg \},
\end{eqnarray}
that is, $\mathcal{T}( \{f_k\}_{k\in \mathbb{Z}})= \{f_{k+1}\}_{k\in \mathbb{Z}}.$
Let us recall that a sequence $\{\frak{M} _k\}_{k \in \mathbb{Z}}$  is  called a  \textit{fusion frame}  if there exist constants
$0 <A \leq B< \infty$ and $v_k > 0$ for $k = 1, 2, \ldots$  so that for every $f \in H$
\begin{eqnarray*} 
A \| f\|^ 2 \leq  \sum_{k \in \mathbb{Z}} v_k^2 \| \pi_{\frak{M} _k} (f) \|^2 \leq B  \|f\|^2, 
\end{eqnarray*}
where $\pi_{\frak{M} _k}$ is the orthogonal projection of $H$ onto $\frak{M} _k$; see  [1, 2, 7]  for more details.
Now, let  $V$ and  $W$  be  Hilbert
spaces and  $\{W_k\}_{k  \in  J} $  be a sequence of closed subspaces of  $W$, where  $J$  is a subset of  $\mathbb{Z}$.
A  sequence  $\{\Lambda_k \}_{ k  \in  J }$ of bounded linear operators  from $V$ into $W_k$  is called a {\it generalized frame  for $ V$ with respect to} $\{W_k\}_{ k \in  J }$  if there are  positive constants $A$ and $B$ such that
\begin{eqnarray*}
A \, \|  f\|^ 2 \leq  \sum_{k\in J} \| \Lambda _k  \, f \|^2 \leq B \, \|f\|^2
\end{eqnarray*}
for all $f \in V$. Note that the fusion frames are a special case of generalized frame; it  is sufficient to set 
\begin{eqnarray*}
 V=W  \,\,\,\,\,\,\,\,\,\,\,\,  \text{and} \,\,\,\,\,\,\,\,\,\,\, \Lambda_k= v_k \pi_{\frak{M} _k}. 
\end{eqnarray*}
A natural way to get an extension of the concept of fusion frame is to replace  the  orthogonal projections on $\frak{M} _k$ by any operator whose range is contained  in $\frak{M} _k$,  $k \in \mathbb{Z}$.

\begin{definition}
Let $\{\frak{M} _k \} _{ k \in \mathbb{Z}} $ be a sequence of closed subspaces of Hilbert space $H$. Let $\{\Theta_k\}_{k \in \mathbb{Z}}$  
be a sequence of operators from  $H$ into $\frak{M} _k$, $k \in \mathbb{Z}$.  The sequence $\{(\frak{M} _k, \Theta_k)\}_{k \in \mathbb{Z}}$ 
 is called a 
$\textit{$g-$fusion frame with respect to}$ $\{\frak{M} _k \} _{ k \in \mathbb{Z}} $  (briefly, $g-$fusion frame) if 
 there exist two constants  $0 < A \leq  B < \infty$ such that for every $f \in H$
\begin{eqnarray}  \nonumber
A  \| f\|^ 2 \leq  \sum_{k \in \mathbb{Z}} \| \Theta _k  (f )\|^2 \leq B  \|f\|^2. 
\end{eqnarray}
In the case, $A=B$, it is called a \emph{tight $g-$fusion frame}. 
\end{definition}
The  bounded operators
$ U : \mathcal{H} \rightarrow  H$ and  $S : H \rightarrow  H$ defined   by
\begin{eqnarray} \nonumber
U \bigg( \{f_k\}_{k \in \mathbb{Z}} \bigg) = \sum_{k\in \mathbb{Z}} \Theta^{*}_{k}  (f_k)\quad \hbox{and}\quad Sf = \sum_{k \in \mathbb{Z} } {\Theta}_{k}  ^{*}   {\Theta}_{k}   f   \nonumber
\end{eqnarray}
are called the $\emph{synthesis operator}$ and  the $\textit{ frame operator}$,  respectively.
A $g-$fusion frame  $\{(\frak{M} _k, \Gamma_k)\}_{k \in \mathbb{Z}}$  is called
a $\emph{dual $g-$fusion frame}$ of   $\{(\frak{M} _k, \Theta_k)\}_{k \in \mathbb{Z}}$   if it satisfies
\begin{eqnarray} \nonumber
f= \sum_{k \in \mathbb{Z}} \Theta_k^{*}  {\Gamma}_k   f,  \,\,\,\,\,\,\,\,\,\,\,\,\,\,\,\,\,\,\,\,\,\,  (f \in H).
\end{eqnarray}
It is easy to see that 
the sequence  $\{(\frak{M} _k, \Theta_k S^{-1})\}_{k \in \mathbb{Z}}$  is  a  dual  $g-$fusion frame of  $\{(\frak{M} _k, \Theta_k)\}_{k \in \mathbb{Z}}$, which is called  $\textit{canonical dual $g-$fusion frame}$.
Note that for any  $f \in H$, we  have
\begin{eqnarray} \nonumber
f = S^{-1}  S  f  =\sum_{k \in \mathbb{Z}} {\widetilde{\Theta_k}}^{*}  \Theta_k  f,  \,\,\,\,\,\,\,\,\, \hbox{where} \,\,\,\,\,
 \widetilde {\Theta}_{k}= \Theta_k   S^{-1}.   \nonumber
\end{eqnarray}

Let $T$ be a linear map on $\hbox{span}\{\frak{M} _k\}_{k \in \mathbb{Z}}$ and $\{\Theta_k\}_{k \in \mathbb{Z}}$ be a sequence in 
$\mathcal{L}(H)$ so that the range of $\Theta_k$ is contained in $\frak{M} _k$, $k \in \mathbb{Z}$. We 
say that  $\{\Theta_k \}_{k \in \mathbb{Z}}$
is $\textit{represented via T}$ if
$\Theta_{k+1}= T  \Theta_k$
for all $ k \in \mathbb{Z}$.
A $g-$fusion frame $\{(\frak{M} _k, \Theta_k)\}_{k \in \mathbb{Z}}$  is said to be represented by a linear  operator $T$ if 
$\{\Theta_k\}_{k \in \mathbb{Z}}$ is represented via $T$. 
One can observe that if  the $g-$fusion frame $\{(\frak{M} _k, \Theta_k)\}_{k \in \mathbb{Z}}$ is represented via a linear operator on $\hbox{span}\{\frak{M} _k\}_{k\in{\mathbb Z}}$, then the same is true for its canonical dual $g-$fusion frame. 
\section{On the  Boundedness of $T$}
In the following,  let $\{\Theta_k\}_{k \in \mathbb{Z}}$ be a sequence in $\mathcal{L}(H)$ such that for every $k \in \mathbb{Z}$, the range of $\Theta_k$ is contained in $\frak{M} _k$. The main result of this section is the following.
\begin{theorem} 
Let  $\Theta_k$  be  self-adjoint and $\Theta_k(\frak{M} _k)= \frak{M} _k$
for all $k \in \mathbb{Z}$.
If  $\{(\frak{M} _k, \Theta_k)\}_{k \in \mathbb{Z}}$ is a 
$g-$fusion frame  represented via a linear operator $T$ on
$\hbox{span} \{\frak{M} _k\}_{ k \in \mathbb{Z}}$, then
 $T$ is bounded and  the kernel of  $U$ is invariant under the  right-shift operator.
Furthermore,
\begin{eqnarray} \nonumber
1 \leq \|T\| \leq \sqrt{B A^{-1}},
\end{eqnarray}
where $A$ and $B$ are  frame  bounds of
 $\{\Theta_k \}_{k \in \mathbb{Z}}$.
\end{theorem}

\begin{proof}
Let $N_U$ denote the kernel of $U$.
If $\{g_k \}_{k \in \mathbb{Z}}$ is a sequence in $N_U^{\perp}$
such that $g_k=0$ for all  but a finite number of $k$, then it
follows  from  Corollary $3.2$ in [14] and Theorem  $5.5.5$ in [3]  that
\begin{eqnarray} \nonumber
\|T U  (\{g_k\}_{k \in \mathbb{Z}} )\|^2 \leq B  \sum_{k \in \mathbb{Z} } \|g_k \|^2 \leq \frac{B}{A}   \|U  ( \{g_k\}_{k \in \mathbb{Z}} ) \|^2.
\end{eqnarray}
This implies that
\begin{eqnarray}
\| T U  (\{g_k\}_{k \in \mathbb{Z}} + \{h_k\}_{k \in \mathbb{Z}}) \|^2 \leq \frac{B}{A}   \| U  (\{g_k \}_{k \in \mathbb{Z}} +\{h_k \}_{k \in \mathbb{Z}}  ) \|^2
\end{eqnarray}
for all  sequence $\{h_k\}_{k \in \mathbb{Z}}$ in $N_U$.
From  (3.1)  and the fact that $\hbox{span}  \{\frak{M} _k\}_{k \in \mathbb{Z}}$  is a subspace of the range of  $U$, we infer that
\begin{eqnarray} \nonumber
\| T  (f)  \|^2 \leq \frac{B}{A}  \| f \|^2
\end{eqnarray}
for all  $f \in \hbox{span}\{\frak{M} _k\}_{k \in \mathbb{Z}} $.
Therefore, $T$ is bounded and $\|T\| \leq \sqrt{B A^{-1}}$. 
Since  $\{(\frak{M} _k, \Theta_k)\}_{k \in \mathbb{Z}}$  is a 
$g-$fusion frame  represented via $T$, there exists $f_0 \in H$ such that  $\sum_{k \in \mathbb{Z}}\|\Theta_k (f_0) \|^2 \not =0$. On the other  hand,
\begin{eqnarray} \nonumber
 \sum_{k \in \mathbb{Z}}    \bigg  \|   \Theta_k  (f_0)   \bigg  \|^2  =  \sum_{k \in \mathbb{Z}}    \bigg  \|     \Theta_{k+1}   (f_0)   \bigg  \|^2  =   \sum_{k \in \mathbb{Z}}    \bigg  \|    T  \Theta_{k}   (f_0)   \bigg  \|^2    \leq \|T\|^2   \sum_{k \in \mathbb{Z}}   \bigg  \|  \Theta_k   (f_0)   \bigg  \|^2.
\end{eqnarray}
This shows that $\|T\| \geq 1$.

To complete the proof, let $\{f_k\}_{k \in \mathbb{Z}} \in N_U$. Then $\sum_{k \in \mathbb{Z}} \Theta_k (f_k ) =0$. Since $T$ is bounded and
$T \Theta_k= \Theta_{k+1}$, we have
\begin{eqnarray*}
\sum_{k \in \mathbb{Z}} \Theta_{k}(f_{k -1}) = \sum_{k \in \mathbb{Z}} \Theta_{k+1} (f_{k}) =0.
\end{eqnarray*}
Hence $U  \mathcal{T} (\{f_k\}_{k \in \mathbb{Z}})=0$.   That  is, $\mathcal{T}(\{f_k\}_{k \in \mathbb{Z}}) \in N_U$.
\end{proof}

As a consequence of Theorem 3.1,  we present the next result.
\begin{corollary}
Let $\{(\frak{M} _k, \Theta_k)\}_{k \in \mathbb{Z}}$   be a $g-$fusion frame represented via a linear operator $T$ on $\hbox{span}\{\frak{M} _k\}_{k \in \mathbb{Z}}$ and $\Theta_k$ be self-adjoint for all $k \in \mathbb{Z}$. If $\Theta_k$ is  a projection on $\frak{M} _k$ for all $k \in \mathbb{Z}$,
then $T$ is bounded.
\end{corollary}
\begin{proof}
Let $\Theta_k$ be  a projection on $\frak{M} _k$ for all $k \in \mathbb{Z}$. Then $\Theta_k(\frak{M} _k)=\frak{M} _k$ for all $k \in \mathbb{Z}$.
This together with Theorem 3.1 proves the result.
\end{proof}
We finish this section with the following result.
\begin{theorem}
Let $\{(\frak{M} _k, \Theta_k)\}_{k \in \mathbb{Z}}$ be a  tight $g-$fusion frame such that  $\Theta_k$ is self-adjoint and $\Theta_k(\frak{M} _k)=\frak{M} _k$
for all $k \in \mathbb{Z}$.
Then  $\{\Theta_k\}_{k \in \mathbb{Z}}$  is not represented via an invertible linear operator $T$ on  $\hbox{span}\{\frak{M} _k\}_{k \in \mathbb{Z}}$.
\end{theorem}
\begin{proof}
Let
$\{ (\frak{M} _k , \Theta_k )\}_{k \in \mathbb{Z}}$ be a $g-$fusion frame  represented via an invertible linear operator $T$ on $\hbox{span} \{\frak{M} _k\}_{k \in \mathbb{Z}}$. Hence
$\Theta_k=T^k\Theta_0$
and so
$\Theta_{-k}=T^{-k}\Theta_0$ for all $k \in \mathbb{Z}$.
Similarly, replacing $T$  by $T^{-1}$ in the proof of Theorem 3.1, we get
\begin{eqnarray} 
1 \leq \|T^{j}\| \leq \sqrt{BA^{-1}},
\end{eqnarray}
where $A$ and $B$ are  frame bounds of $\{ (\frak{M} _k , \Theta_k )\}_{k \in \mathbb{Z}}$ and $j\in\{-1, 1\}$.
Since $\{ (\frak{M} _k , \Theta_k )\}_{k \in \mathbb{Z}}$ is a tight $g-$fusion frame,
then  $A=B$. So by (3.2)
\begin{eqnarray} \nonumber
 \|T\|=\|T^{-1}\|=1.
\end{eqnarray}
Therefore  for every $f \in H$,
\begin{eqnarray} \nonumber
\|f\|= \|T ^{-1}  T f\| \leq \|Tf\| \leq \|f\|.
\end{eqnarray}
This shows that  $T$ is an isometry.
Thus  for every $f \in H$ and $k \in \mathbb{Z}$
\begin{eqnarray} \nonumber
\|\Theta_k  f \|= \|T^k  \Theta_0  f \| = \| \Theta_0 f \|.
\end{eqnarray}
Hence
\begin{eqnarray} \nonumber
\sum_{k\in {\mathbb Z}} \|\Theta_0(f)\|^2 =\sum_{k\in {\mathbb Z}}\|\Theta_k(f)\|^2 \leq B \|f\|^2.
\end{eqnarray}
It follows that
$\sum_{k\in {\mathbb Z}}\|\Theta_0(f)\|^2$
is a convergent series. So
\begin{eqnarray} \nonumber
\|\Theta_k(f)\|=\|\Theta_0(f)\|=0
\end{eqnarray}
for all $f\in H$ and $k \in \mathbb{Z}$.
Therefore  for every $f\in H$ we have
$$
f=\sum_{k\in {\mathbb Z}}\Theta_k^*\Gamma_k f=0,$$
where $\{(\frak{M} _k, \Gamma_k) \}_{k \in \mathbb{Z}}$ is a dual
$g-$fusion frame of $\{(\frak{M} _k, \Theta_k)\}_{k \in \mathbb{Z}}$, a contradiction.
\end{proof}
\section{Linear Independence and Stability}

We start with a result  that plays an important role in the sequel.
\begin{proposition}
Let $\{\Theta_k \}_{k \in \mathbb{Z}} $  be a sequence in $\mathcal{L}(H)$ with the range of $\Theta_k$ is contained in $\frak{M} _k$ for all $k \in \mathbb{Z}$. If 
$\{\Theta_k \}_{k \in \mathbb{Z}} $ is represented via a  linear operator $T$ on
 $ \hbox{span} \{\frak{M} _k \}_{k \in \mathbb{Z}} $
and   $ \hbox{span} \{\Theta_k \}_{k \in \mathbb{Z}} $
is   infinite dimensional, then 
$ \{\Theta_k \}_{k \in \mathbb{Z}} $
 is linearly independent and infinite.
\end{proposition}
\begin{proof}
Suppose that
$\{\Theta_k \}_{k \in \mathbb{Z}} $
is linearly dependent. Hence there exist constants
$c_M, \ldots, c_N$ such that
\begin{eqnarray} \nonumber
\sum_{k=M}^{N} c_k \Theta_k =0
\end{eqnarray}
and  $ c_{k_0} \not = 0$, for some $ M \leq k_0 \leq N$. Set
\begin{eqnarray} \nonumber
\alpha= \min \bigg\{ k : \,\,\,  M \leq k \leq N, c_k \not =0 \bigg \} \,\,\,\, \hbox{and}  \,\,\,\,\,
\beta= \max \bigg\{ k : \,\,\,  M \leq k \leq N, c_k \not =0 \bigg \}.
\end{eqnarray}
Then
\begin{eqnarray} \nonumber
\Theta_{\beta} = \sum_{k =\alpha}^{\beta-1} d_k  \Theta_k  \,\,\,\,\,\,\,\,\, \,\,\,\,\,\,\,\,\,   \hbox{and}  \,\,\,\,\,\,\,\,\, \,\,\,\,\,\,\,\,\,
\Theta_{\alpha} = \sum_{k =\alpha+1}^{\beta} d_k^{\prime} \Theta_k
\end{eqnarray}
 for some constants
$d_k$ and $d_k^{\prime}$, $k= \alpha, \ldots, \beta$.
Thus for every $i \in \mathbb{N}$, we have
\begin{eqnarray} \nonumber
T^i  \Theta_{\beta} = \sum_{k=\alpha} ^{\beta -1} d_k  T^i  \Theta_k
=\sum_{k=\alpha} ^{\beta -1} d_k  \Theta_{i+k }
= \sum_{k=\alpha+i} ^{\beta+ i -1} d_{k-i} \Theta_{k }.
\end{eqnarray}
Similarly,
$T^{-i}  \Theta_{\beta} =  \sum_{k=\alpha-i} ^{\beta-1-i } d_{k-i}^{'}   \Theta_{k }$.
These facts together with
$\Theta_k= T^k  \Theta_0$ show that
\begin{eqnarray} \nonumber
{\cal V} :=  \hbox{span} \big \{ \Theta_k \big \}_{ \alpha, \ldots, \beta}
\end{eqnarray}
is invariant under $T$ and $T^{-1}$. Also, we have
\begin{eqnarray} \nonumber
 \hbox{span} \big\{ \Theta_k  \big\}_{ k \in \mathbb{Z} }={\cal V}.
\end{eqnarray}
Therefore,
$ \hbox{span} \big \{ \Theta_k \big\}_{k \in \mathbb{Z}}$
is  finite dimensional, a contradiction.
\end{proof}

As an immediate consequence of Proposition 4.1, we have the following results.
\begin{theorem} 
Let  there exist $j_0 \in \mathbb{Z}$ such that $\Theta_{j_0}$ is surjective  and for every $k \in \mathbb{Z}$, $\Theta_k$ is self-adjoint and $\Theta_k(\frak{M} _k)=\frak{M} _k$.
Let $\{ (\frak{M} _k , \Theta_k )\}_{k \in \mathbb{Z}}$ be a $g-$fusion frame   and  $T$ be a linear operator on $\hbox{span}\{\frak{M} _k\}_{k\in{\mathbb Z}}$ such that it has an extension to a  bounded linear operator $\tilde{T}: H\rightarrow H$ with $\tilde{T}
\Theta_k^*=\Theta_{k+1}^*$.
If $\{\Theta_k \}_{k \in \mathbb{Z}}$ is represented via  $T$,
then $\{\Theta_k \}_{k \in \mathbb{Z}}$ is linearly dependent.
\end{theorem}
\begin{proof}
First, let   $\{(\frak{M} _k, \Gamma_k )\}_{k \in \mathbb{Z}}$  be the dual   $g-$fusion frame of 
$\{ (\frak{M} _k , \Theta_k )\}_{k \in \mathbb{Z}}$. Thus  every $f \in H$ has the decomposition
\begin{eqnarray} \nonumber
f = \sum_{k \in \mathbb{Z}}  \Theta_{k}^{*}    \Gamma_k  (f).
\end{eqnarray}
So $\tilde{T}(f)= \sum_{k \in \mathbb{Z}}  \Theta_{k+1}^{*}    \Gamma_k   (f)$. This implies that
\begin{eqnarray*} \nonumber
T(\Theta_j (f) ) =  \sum_{k \in \mathbb{Z}}  \Theta_{k+1}^{*}    \Gamma_k    \Theta_j  (f)
\end{eqnarray*}
and then for every $f \in H$ we have
\begin{eqnarray} 
\Theta_{j+1}  (f)= \sum_{k \in \mathbb{Z}}  \Theta_{k+1}^{*}    \Gamma_k    \Theta_j  (f).
\end{eqnarray}
Choose $\ell, \ell^{'} \in \mathbb{Z}$ such that  $\ell \not = \ell^{'}$ and $ j_{0}  \not \in  \{\ell -1, \ell,\ell^{'}-1, \ell^{'} \}$.
Let $\{\tilde{\Theta}_k\}_{k \in \mathbb{Z}}$ denote the sequence consisting of
the same elements as $\{\Theta_k\}_{k \in \mathbb{Z}}$, but with $\Theta_\ell$ and $\Theta_{\ell^\prime}$ interchanged.
It is clear that
\begin{eqnarray} \nonumber
\widetilde{\Gamma}_k =
\left\{
  \begin{array}{ll}
    \Gamma_k, \,\,\,\,\,\,\,\,\,\,\,\,\,\,    k \not = \ell, \ell^{'},  \\
    \Gamma_{\ell^{'}},  \,\,\,\,\,\,\,\,\,\,\,\,\,\,     k=\ell, \\
    \Gamma_{\ell},  \,\,\,\,\,\,\,\,\,\,\,\,\,\,     k=\ell^{'},
  \end{array}
\right.
\end{eqnarray}
is the canonical dual $g-$fusion frame $\{(\frak{M} _k, \widetilde{\Theta}_k ) \}_{ k \in \mathbb{Z}}$.
 From (4.1) we obtain
\begin{eqnarray} 
\Theta_{j_{0}+1} &=&  \sum_{k \in \mathbb{Z}} \Theta_{k+1}  {\Theta_k  S^{-1}}   \Theta_{j_{0}}  \nonumber \\
&=&  \sum_{k \not \in  \{\ell -1, \ell, \ell^{'}-1, \ell{'}\} }\Theta_{k+1}   {\Theta_k  S^{-1}}   \Theta_{j_{0}}     \\
&+&  \Theta_{\ell}  \Theta_{\ell-1}  S^{-1}   \Theta_{j_{0}}
+  \Theta_{\ell+1}   \Theta _{\ell}   S^{-1}    \Theta_{j_{0}}  \nonumber \\
&+ & \Theta_{\ell^{'}}   \Theta_{\ell^{'}-1}   S^{-1}    \Theta_{j_{0}}
+\Theta_{\ell^{'}+1}  \Theta_{\ell^{'}}   S^{-1}   \Theta_{j_{0}}. \nonumber
\end{eqnarray}
Then we have
\begin{eqnarray} 
\Theta_{j_{0}+1} &=&  \widetilde{\Theta}_{j_{0}+1} \nonumber \\
&=&  \sum_{k \in \mathbb{Z}}  \widetilde{\Theta}_{k+1}  \widetilde{\Gamma}_k    \widetilde{\Theta }_{j_{0}}  \nonumber \\
& =&   \sum_{k \not \in  \{\ell -1, \ell, \ell^{'}-1, \ell{'}\} }\Theta_{k+1}  \Theta_k   S^{-1}  \Theta_{j_{0}} \nonumber \\
&+&  \Theta_{\ell ^{'}}   \Theta_{\ell-1}  S^{-1}   \Theta_{j_{0}}
+  \Theta_{\ell+1} \Theta_{\ell^{'}}  S^{-1}   \Theta_{j_{0}} \nonumber \\
&+& \Theta_{\ell}   \Theta_{\ell^{'}-1}  S^{-1}  \Theta_{j_{0}}
+ \Theta_{\ell^{'}+1}   \Theta_{\ell}  S^{-1}  \Theta_{j_{0}}. \nonumber
\end{eqnarray}
Using   (4.2), we obtain
\begin{eqnarray} \nonumber
 \Theta_{\ell}  \Theta_{\ell-1}  S^{-1}   \Theta_{j_{0}}
&+&  \Theta_{\ell+1}   \Theta _{\ell}   S^{-1}    \Theta_{j_{0}}
+  \Theta_{\ell^{'}}   \Theta_{\ell^{'}-1}   S^{-1}    \Theta_{j_{0}} \nonumber \\
&+& \Theta_{\ell^{'}+1}   \Theta_{\ell^{'}}   S^{-1}   \Theta_{j_{0}}
=   \Theta_{\ell ^{'}}   \Theta_{\ell-1}  S^{-1}   \Theta_{j_{0}}
\nonumber \\
&+&  \Theta_{\ell+1}   \Theta_{\ell^{'}}  S^{-1}   \Theta_{j_{0}}
+ \Theta_{\ell}   \Theta_{\ell^{'}-1}  S^{-1}   \Theta_{j_{0}}  \nonumber \\
&+& \Theta_{\ell^{'}+1}   \Theta_{\ell}   S^{-1}  \Theta_{j_{0}}. \nonumber
\end{eqnarray}
Without loss of generality, we may assume that $\ell >\ell^{'}$ and $\ell^{'}= \ell+N$,  for some $N \in \mathbb{N}$.
So
\begin{eqnarray} \nonumber
 \Theta_{\ell}   \Theta_{\ell-1}   S^{-1}   \Theta_{j_{0}}
&+&  \Theta_{\ell+1}  \Theta_{\ell} S^{-1}  \Theta_{j_{0}}
+ \Theta_{\ell+N}  \Theta_{\ell+(N-1)}  S^{-1}  \Theta_{j_{0}} \nonumber \\
&+& \Theta_{\ell+(N+1)}  \Theta_{\ell+N}  S^{-1}   \Theta_{j_{0}}
- \Theta_{\ell+N} \Theta_{\ell-1}   S^{-1}   \Theta_{j_{0}} \nonumber \\
&-&  \Theta_{\ell+1}  \Theta_{\ell+N}   S^{-1}   \Theta_{j_{0}}
- \Theta_{\ell}    \Theta_{\ell+(N-1)}  S^{-1} \Theta_{j_{0}} \nonumber \\
&-& \Theta_{\ell+(N+1)}   \Theta_{\ell}   S^{-1}  \Theta_{j_{0}}=0.   \nonumber
\end{eqnarray}
Hence
\begin{eqnarray} \nonumber
  \Theta_{2  \ell-1}   S^{-1}   \Theta_{j_{0}}
&+&  \Theta_{2   \ell+1}  S^{-1}  \Theta_{j_{0}}
+ \Theta_{2 (\ell + N) -1}   S^{-1}  \Theta_{j_{0}}  \nonumber \\
&+& \Theta_{2 ( \ell+N)+1}  S^{-1}   \Theta_{j_{0}}
-2    \Theta_{2  \ell+(N-1)}  S^{-1}  \Theta_{j_{0}} \nonumber \\
&-& 2   \Theta_{2  \ell+(N+1)}  S^{-1}   \Theta_{j_{0}} = 0.   \nonumber
\end{eqnarray}
Since  $\Theta_{j_0}$ and $S^{-1}$  are surjective, we have
\begin{eqnarray} \nonumber
  \Theta_{2  \ell-1}  &+& \Theta_{2   \ell+1} + \Theta_{2  (\ell+N)-1} \nonumber \\
&+&   \Theta_{2 ( \ell+N)+1}
-  2  \Theta_{2 \ell+(N-1)}  \nonumber \\
&- & 2  \Theta_{2  \ell+(N+1)} =0,   \nonumber
\end{eqnarray}
so $\{\Theta_k\}_{k \in \mathbb{Z}}$ is linearly dependent. 
\end{proof}
\begin{theorem} 
Let
$\{ (\frak{M} _k , \Theta_k )\}_{k \in \mathbb{Z}}$  be a $g-$fusion frame  represented via a linear operator $T$ on $\hbox{span}\{\frak{M} _k\}_{k\in {\mathbb Z}}$   and
 $\{\widehat{\Theta}_k \}_{k \in \mathbb{Z}}$  be a sequence in ${\cal L}(H)$ such that $\widehat{\Theta}_k(H) \subseteq \frak{M} _k$ for all $k \in \mathbb{Z}$.  If there exist $\alpha, \beta \in [0,1)$ such that
\begin{eqnarray}
\bigg \| \sum_{k=1}^ n  c_ k\big(\Theta_k - \widehat{\Theta}_k  \big) (f) \bigg \| \leq \alpha
\bigg \| \sum_{k=1}^ n  c_ k  \Theta_k  (f) \bigg \| + \beta
\bigg \| \sum_{k=1}^ n  c_ k  \widehat{\Theta}_k (f) \bigg \|
\end{eqnarray}
for all  $f \in H$ and finite  complex sequences $\{c_k\}_{k=1}^n$, then the following statements hold.
\begin{itemize}
\item[$\emph{(i)}$]
$\{(\frak{M} _k, \widehat{\Theta}_k )\}_{k \in \mathbb{Z}}$  is  a $g-$fusion frame   with frame bounds
$(\frac{1- \alpha}{1+ \beta} )^2  A$ and $(\frac{1+ \alpha}{1- \beta} )^2 B$, where $A$ and $B$ are  frame bounds  of $\{ (\frak{M} _k , \Theta_k )\}_{k \in \mathbb{Z}}$ . 
\item[$\emph{(ii)}$]
If $\hbox{span}\{\Theta_k\}_{k \in \mathbb{Z}}$ is infinite dimensional and $\{\widehat{\Theta}_k \}_{k \in \mathbb{Z}}$  is  an  infinite set, then
$\{\widehat{\Theta}_k \}_{k \in \mathbb{Z}}$ is linearly independent.
\end{itemize}
\end{theorem}
\begin{proof}
Choose $i \in \mathbb{Z}$ and set $c_k = \varepsilon_{ik}$  in  (4.3), where $ \varepsilon_{ik}$ is the Kronecker delta. Hence for every $f \in H$
\begin{eqnarray} \nonumber
\| \Theta_i (f) - \widehat{\Theta}_i (f) \| \leq \alpha  \| \Theta_i (f)  \|  + \beta  \| \widehat{\Theta}_i (f) \|
\end{eqnarray}
and then
\begin{eqnarray} \nonumber
(1 - \beta)   \| \widehat{\Theta}_i (f) \| \leq (1 + \alpha)  \| \Theta_i (f)  \|.
\end{eqnarray}
Therefore
\begin{eqnarray} \nonumber
 \| \widehat{\Theta}_i (f) \|^2 \leq  \bigg(\frac{1 + \alpha}{1 - \beta} \bigg)^2   \| \Theta_i (f)  \|^2.
\end{eqnarray}
This shows that
\begin{eqnarray} \nonumber
\sum_{ i \in \mathbb{Z}}  \| \widehat{\Theta}_i (f) \|^2 \leq  \bigg(\frac{1 + \alpha}{1 - \beta} \bigg)^2   \sum_{ i \in \mathbb{Z}}   \| \Theta_i (f)  \|^2 \leq  \bigg(\frac{1 + \alpha}{1 - \beta} \bigg)^2   B  \|f\|^2,
\end{eqnarray}
where $B$ is  the upper frame bound of $\{ (\frak{M} _k , \Theta_k )\}_{k \in \mathbb{Z}}$.
Similarly, one can show that
$ \| \widehat{\Theta}_i (f) \|^2 \geq  \bigg(\frac{1 - \alpha}{1 + \beta} \bigg)^2  \| \Theta_i f\|^2$
 for all $i \in \mathbb{Z}$ and $f \in H$.
So
\begin{eqnarray} \nonumber
\sum_{ i \in \mathbb{Z}}  \| \widehat{\Theta}_i (f) \|^2 \geq  \bigg(\frac{1 - \alpha}{1 + \beta} \bigg)^2  A   \|f\|^2,
\end{eqnarray}
 where
$A$ is  the lower  frame  bound of $\{ (\frak{M} _k , \Theta_k )\}_{k \in \mathbb{Z}}$. Therefore,
$\{ (\frak{M} _k , \widehat{\Theta}_k )\}_{k \in \mathbb{Z}}$  is  a $g-$fusion frame  with  frame  bounds
\begin{eqnarray} \nonumber
 \bigg(\frac{1 - \alpha}{1 + \beta} \bigg)^2  A \,\,\,\,\,\,\,\, \hbox{and} \,\,\,\,\,\,\,\, \bigg(\frac{1 + \alpha}{1 - \beta} \bigg)^2  B.
\end{eqnarray}
That is, $\hbox{(i)}$ holds.

Let $\{c_k \}_{k=1}^n$ be a finite sequence in $\mathbb{C}$ such that  $\sum_{k=1}^n c_k  \widehat{\Theta}_k =0$. Since
$\alpha \in [0,1)$,  by  (4.3)
\begin{eqnarray} \nonumber
\sum_{k=1}^n c_k \Theta_k =0.
\end{eqnarray}
By Proposition 4.1, $c_k =0$ for all $k=1, \ldots, n$.
So $\{\widehat{\Theta}_k \}_{k \in \mathbb{Z}}$ is linearly independent.
\end{proof}
\section{Some Examples}

Let ${\cal G}$ be a locally compact Hausdorff group with a left Haar measure $\lambda$ and identity element $e$.
We denote by $M({\cal G})$ the space of complex regular Borel measures on $\mathbb{C}$. For every $\mu\in M({\cal G})$ and $f\in L^2({\cal G}):=L^2({\cal G}, \lambda)$ we define $\mu\ast f \in L^2({\cal G})$ by
\begin{eqnarray} \nonumber
\mu\ast f(x)=\int_{\cal G} f(y^{-1}x) d  \mu(y)
\end{eqnarray}
for all $x\in {\cal G}$.
We recall  that for every $x\in {\cal G}$ and $f\in L^2({\cal G})$
\begin{eqnarray*}
\|\delta_x\ast f\|_2=\|f\|_2,
\end{eqnarray*}
where $\delta_x$ is the  Dirac measure at $x$  [9].
\begin{example} 
Let ${\cal G}$ be an infinite locally compact Hausdorff group. Choose a sequence $\{x_k\}_{k \in \mathbb{Z}}$ in ${\cal G}$ such that $x_0=e$,
$x_k=x_{-k}^{-1}$
 for all $k<0$ and $x_k \not = x_{\ell}$, for $k \not = \ell$.
 Define $\Theta_k: L^2({\cal G}) \rightarrow L^2({\cal G})$ by
\begin{eqnarray*}
\Theta_k(f)= (1/2)^{|k|} \delta_{x_k}\ast f
\end{eqnarray*}
for all $k\in{\mathbb Z}$.
Let $\sum_{k=M}^N c_k  \Theta_k =0$ for some $M, N  \in \mathbb{Z}$ and $c_k \in \mathbb{C}$.
 Since ${\cal G}$ is Hausdorff, there is an open set $O$ with $x_M^{-1} \in O$  and $x_k^{-1} \not \in O$ for $k=M+1, \ldots, N$.
By  Urysohn's  Lemma, there exists $f \in L^2({\cal G})$ such that $f (x_M^{-1})=1$ and $f(x_k^{-1}) =0$ for all $k=M+1, \ldots, N$. Thus
\begin{eqnarray} \nonumber
0= \sum_{k =M}^{N} c_k  \Theta_k (f)  (e) =  \sum_{k =M}^{N} c_k   f(x_k^{-1}) = c_M.
\end{eqnarray}
Similarly, $c_k =0$ for $k=M+1, \ldots, N$.
 Hence $\{ \Theta_k\}_{k \in \mathbb{Z}}$ is  linearly independent. Note that $\{ (\frak{M} _k , \Theta_k )\}_{k \in \mathbb{Z}}$ is a tight $g-$fusion frame for $L^2({\cal G})$. 
If $\{ \Theta_k\}_{k \in \mathbb{Z}}$ is represented  via  a linear operator $T$ on $\hbox{span}\{\frak{M} _k\}_{k \in \mathbb{Z}}$, then  $\|T\|_2= 1/2$.
Indeed, for every $f \in L^2 ({\cal G})$ we have $\|\Theta_k(f)\|_2=(1/2)^{|k|}  \|f\|_2$ and
\begin{eqnarray} \nonumber
\|T  (f) \|_2  &=& \| T  ( \delta_{x_1} \ast   \delta_{{x}^{-1}_1} \ast   f) \|_2  \nonumber \\
&=& 2\;\|T  \Theta_1  ( \delta_{{x}^{-1}_1} \ast  f) \|_2 \nonumber \\
& =&2\;\| \Theta_2   ( \delta_{{x}^{-1}_1} \ast   f) \|_2 \nonumber \\
& =& (1/2)\; \|  \delta_{x_2} \ast   \delta_{x^{-1}_1} \ast   f \|_2  \nonumber \\
&=& (1/2)\; \|f\|_2.  \nonumber
\end{eqnarray}
\end{example}

The next example  shows that the requirement that $\Theta_k$ be self-adjoint on $\frak{M} _k$ in Theorems 3.1 and 3.3 can not be omitted.
\begin{example}
 Let ${\cal G}$ be a discrete group. Then $L^2({\cal G})$ is isometrically isomorphic with $\ell^2({\mathbb Z})$.
Let $\{\Theta_k\}_{k \in \mathbb{Z}}$  be as in Example 5.1. As was noted,
$\{ (\frak{M} _k , \Theta_k )\}_{k \in \mathbb{Z}}$  is a tight $g-$fusion frame represented by a bounded linear operator $T$ and
   $\|T\|= 1/2$. It is easy to see that
\begin{eqnarray} \nonumber
\Theta_k^{\ast}(f)= (1/2)^{|k|} \delta_{x^{-1}} \ast f
\end{eqnarray}
for all $f \in L^2({\cal G})$.
Hence $\Theta_k$ is not self-adjoin on $\frak{M} _k$.
\end{example}
\begin{example}
Let ${\cal G}$ be a locally compact Hausdorff group and $\{\frak{M} _k \}_{k \in \mathbb{Z}}$ be a sequence of   closed subspaces of $L^2({\cal G})$.
We define $\Phi: M({\cal G}) \rightarrow \mathcal{L}(L^2 ({\cal G}))$ by
\begin{eqnarray*}
 \Phi (\mu) (f)= \mu\ast f.
\end{eqnarray*}
Choose $x_0 \in {\cal G}$ and set $T= \Phi (\delta_{x_0})$.
We will  show that there is no  $g-$fusion frame for $L^2({\cal G})$  that can be represented via $T$. To do this, suppose that
$\{ (\frak{M} _k , \Theta_k )\}_{k \in \mathbb{Z}}$ is a  $g-$fusion frame  for
  $L^2({\cal G})$ such  that
 $\{\Theta_k \}_{ k \in \mathbb{Z}}$ is represented via $T$. Then for every $f \in L^2({\cal G})$
\begin{eqnarray} \nonumber
\Theta_k (f) = T^k  \Theta_0 (f) = \delta_{x_0^k} * \Theta_0 (f).
\end{eqnarray}
Since $\lambda$ is invariant under left translations, then
\begin{eqnarray} \nonumber
\|\Theta_k   (f)\|_2^2 &=& \int_{\cal G} |  \Theta_k  (f)   |^2  (x)  d  \lambda  (x) \\ \nonumber
&=& \int_{\cal G} |   \delta_{x_0^k} * \Theta_0   (f)  |^2  (x)  d  \lambda  (x) \nonumber \\
&=& \int_{\cal G} \bigg|\int_{\cal G} \Theta_0(f) (y^{-1} x)|^2 d \delta_{x_k}(y) \bigg| d \lambda (x) \nonumber \\
&=& \int_{\cal G}   |  \Theta_0   (f)   (x_0 ^{- k}  x) |^2  d  \lambda  (x)   \nonumber \\
&=& \int_{\cal G} | \Theta_0   (f)  (x)   |^2   d   \lambda  (x) \nonumber \\
&=& \| \Theta_0   (f) \|_2^2.  \nonumber
\end{eqnarray}
Therefore,
\begin{eqnarray} \nonumber
\sum_{ k \in \mathbb{Z}} \| \Theta_k   (f)\|_2^2 = \sum_{ k \in \mathbb{Z}} \| \Theta_0   (f)\|_2^2 = \infty.
\end{eqnarray}
That is, $\{ (\frak{M} _k , \Theta_k )\}_{k \in \mathbb{Z}}$  is not a  $g-$fusion frame for $L^2({\cal G})$. 
\end{example}
\section{Acknowledgments}
The authors thank Professor Robert Brown for assistance with the exposition in this paper. 
The first author spent his sabbatical  at the University of California, Los Angeles (UCLA) during the writing process of this joint work.
\section{Compliance with ethical standards}
\textbf{Conflict of interest:} All authors declare that they have no conflict of
interest.\\
\textbf{Data Availability Statement:} No datasets were generated or analysed during the current study. 
\vspace{0.5 cm}

\end{document}